\documentclass[reqno,12pt]{article}
\usepackage{amsmath}
\usepackage{amssymb}
\usepackage{amsthm}
\usepackage[frenchb]{babel}

\usepackage{epsfig}
 
 \newtheorem{theorem}{Th\'eor\`eme}[section]

\newtheorem{dem}[theorem]{D\'emonstration}
  
\begin{document}

\title {OBSTRUCTION A LA TRANSVERSALITE DE LAGRANGIENS }

\author{Max KAROUBI\,\footnote{  Math\'ematiques, UMR $ 7586$ du CNRS, case $7012$, Universit\'e Paris $7$, $2$, Place Jussieu, $75251$ Paris cedex $05$, France, e.mail : karoubi@math.jussieu.fr.}\and Maria Luiza LAPA de SOUZA \footnote{}\,\footnote { Matem\'atica, Universidade Federal da Bahia, Campus de Ondina, $40170.100$, Ondina, SSA/Bahia, Brasil, e.mail : lapa@noos.fr.}}

\date{}
\maketitle
 
\bigskip
\noindent
{\bf R\'esum\'e.} Dans cette Note, nous donnons une condition n\'ecessaire et suffisante pour que
plusieurs lagrangiens dans un fibr\'e vectoriel symplectique puissent \^etre d\'eform\'es
stablement en lagrangiens transversaux. Dans le cas de trois lagrangiens, nous
montrons que le groupe de Grothendieck associ\'e \`a ce probl\`eme s'identifie \`a un groupe
de K-th\'eorie hermitienne.

\bigskip

OBSTRUCTION TO LAGRANGIAN TRANSVERSALITY

\medskip
\noindent
{\bf Abstract.}
In this Note, we give a necessary and sufficient condition for lagrangians in a
symplectic vector bundle to be deformed stably into transversal lagrangians. In the case
of three lagrangians, we show that the associated Grothendieck group can be identified
with a hermitian K-theory group.

\medskip
\noindent
{\bf 1.} Soit $E$ un fibr\'e vectoriel r\'eel de base $X$ muni d'une forme symplectique non
d\'eg\'en\'er\'ee (dite aussi $(-1)$-hermitienne). Un sous-fibr\'e $L$ est dit lagrangien si son
orthogonal est \'egal \`a L. Deux lagrangiens $L_1$ et $L_2$ sont transversaux si leur intersection
$L_1 \cap L_2$ est r\'eduite \`a $0$. Ils sont homotopiquement transversaux s'il existe deux
familles continues de lagrangiens $L_1$(t) et $L_2$(t) telles que $L_1$(0) = $L_1$, $L_2$(0) = $L_2$, les
lagrangiens $L_1$(1) et $L_2$(1) \'etant transversaux. Enfin, $L_1$ et $L_2$ sont dits stablement et
homotop"quement transversaux si $L_1 \oplus L$ et $L_2 \oplus L$ sont homotopiquement 
transversaux dans le fibr\'e $E \oplus H(L)$, o\`u $H(L) = L \oplus L^*$ d\'esigne le fibr\'e hyperbolique
associ\'e \`a $L$.

\medskip
\noindent
{\bf 2.} Consid\'erons le mono\"\i de (pour la somme directe) des triplets $(E, L_1, L_2)$, o\`u $E $ est un fibr\'e symplectique avec
deux lagrangiens $L_1$ et $L_2$. Un tel triplet est \'el\'ementaire si $L_1$ = $L_2$. Deux triplets $\sigma $ et $\sigma '$
sont \'equivalents s'il existe $\tau $ et $\tau ' $ \'el\'ementaires tels que $\sigma + \tau $ soit isomorphe \`a $\sigma ' + \tau '$
en un sens \'evident. L'ensemble des classes d'\'equivalence forme un groupe $U(X)$ \'etudi\'e
dans \cite{K1} et  \cite{K2} dans un contexte plus g\'en\'eral. On note $ d(E, L_1, L_2)$ la classe du triplet
$(E, L_1, L_2)$ dans $U(X)$. Ce groupe $U(X)$ s'ins\`ere dans une suite exacte entre des groupes de K-th\'eorie r\'eelle et
complexe, o\`u la derni\`ere fl\`eche est induite par la complexification des fibr\'es
$$K_R
^{-1}(X) \rightarrow K_C
^{-1}(X) \rightarrow  U(X)\rightarrow  K_R(X) \rightarrow  K_C(X)$$
En fait, on d\'emontre dans \cite{K1} que $U(X)$ s'identifie \`a un certain groupe $ V^1(X)$ obtenu \`a
partir de triplets $(F, g_1, g_2),$ o\`u $F $ est un fibr\'e r\'eel (sur une $7-$i\`eme-suspension de $X^+$) et o\`u $g_1$
et $g_2$ sont deux formes quadratiques non d\'eg\'en\'er\'ees sur $F$. Plus pr\'ecis\'ement, $U(X) \cong V^1(X)$ est aussi isomorphe \`a $K_R
^1 (X) $ qui est la $K$-th\'eorie associ\'ee \`a la cat\'egorie des
fibr\'es vectoriels r\'eels sur la $7$-i\`eme suspension topologique de $X^+$ (ceci est une version r\'eelle de la
p\'eriodicit\'e de Bott).

\medskip
\noindent
Dans sa th\`ese \cite {LdeS}, le deuxi\`eme auteur a calcul\'e l'obstruction \`a la transversalit\'e
homotopique stable des lagrangiens $L_1$ et $L_2$ , probl\`eme initialement pos\'e au d\'ebut. 
Elle se situe dans le conoyau de la fl\`eche $K_R(X) \rightarrow U(X)$ qui associe \`a un fibr\'e $L$
la classe du triplet $ (H(L), L, L^* )$. En raison de la p\'eriodicit\'e de Bott de nouveau, on
montre que cette fl\`eche est nulle. L'obstruction cherch\'ee appartient donc simplement au
groupe $U(X)\cong K_R
^1 (X)$. Comme il est montr\'e dans \cite {LdeS}, l'invariant de Maslov classique
associ\'e \`a deux lagrangiens se d\'eduit de cet invariant.

\medskip
\noindent
{\bf 3.} Remarque. Si $Q$ est le corps des nombres rationnels, le groupe $K_R
^1 (X) \otimes Q$ s'identifie \`a la somme des groupes de cohomologie de Cech $H^{4n+1}(X ; Q)$.

\medskip
\noindent
{\bf 4.} Examinons maintenant la situation de trois lagrangiens $L_1$, $L_2$ et $ L_3$ avec une
extension \'evidente des d\'efinitions de transversalit\'e. 

\medskip
\noindent
\begin{theorem} Pour que trois lagrangiens $L_1$, $L_2$ et $L_3$ soient homotopiquement et
stablement transversaux, il faut et il suffit que $d(E, L_1, L_2) = d(E, L_1, L_3) = 0$  dans le
groupe $U(X)$ (ce qui implique aussi $d(E, L_2, L_3) = 0$).
\end{theorem}

\noindent
\begin{dem} La condition \'etant \'evidemment n\'eces\-saire, d\'emontrons qu'elle est
suffisante. Puisque $d(E, L_1, L_2) = 0$, on peut supposer (stablement) que $E = H(L)$ avec
$L_2 = L$ et $L_1 = L^*$ qui est isomorphe \`a $L$ en tant que fibr\'e r\'eel par le choix d'une
m\'etrique sur $L$. Par ailleurs, puisque $d(E, L_2, L_3) = 0$, on peut aussi supposer que $L_3 $ est
transverse \`a $L_2$ et qu'il est d\'efini par le graphe d'un morphisme
$g : L = L_1 \rightarrow L^* = L_2.$
Puisque $L_3$ est un lagrangien dans un fibr\'e symplectique, on a n\'ecessairement $g^* = g$.
Soit maintenant $h$ une m\'etrique sur $L$ d\'efinissant donc un isomorphisme auto-adjoint
entre $L$ et $ L^*$. Alors le graphe de l'endomorphisme auto-adjoint $g(t) = (1-t) g + t h$
d\'efinit une homotopie entre $L_3$ et un lagrangien transversal \`a la fois \`a $L_1$ et $L_2$ (pour $ t =
1$).
\end{dem}

\medskip
\noindent
{\bf 5.} La m\^eme m\'ethode permet de d\'emontrer le th\'eor\`eme plus g\'en\'eral suivant.

\medskip
\noindent
\begin{theorem} Soit E un fibr\'e symplectique et soient $L_1$, ..., $L_m$ des lagrangiens
avec $m \geq 2$. Pour que ceux-ci soient homotopiquement et stablement transversaux, il
faut et il suffit que dans le groupe $U(X)$ on ait  $d(E, L_1, L_n) = 0$ pour tout  $n$ (ce qui
implique $d(E, L_i, L_j) = 0$ pour tout couple $(i, j))$.
\end{theorem}

\noindent 
\begin{dem} Le raisonnement pr\'ec\'edent montre que les $L_n$ , $n\geq 3$, peuvent \^etre
associ\'es (\`a homotopie pr\`es) \`a des m\'etriques $h_n$ sur $L = L_1$. Par homoth\'etie, on peut
m\^eme supposer que $ h_n = n h$, o\`u $h$ est une m\'etrique fixe sur $L$. On en d\'eduit que
$h_n - h_p = (n - p) h$ est inversible pour $n \not= p$ et et que donc $L_n$ est transversal \`a $L_p$.
\end{dem}

\medskip
\noindent
{\bf 6.} Nous allons voir que la situation est l\'eg\`erement diff\'erente pour des fibr\'es E munis
d'une forme bilin\'eaire sym\'etrique non d\'eg\'en\'er\'ee (ou encore $(+1)$-hermitienne). Dans ce
cas, le groupe $U(X)$ construit pr\'ec\'edemment avec des couples de lagrangiens est
isomorphe \`a $K_R
^{-1}(X)$ d'apr\`es \cite{K1}. L'obstruction pour la transversalit\'e homotopique
stable de deux lagrangiens appartient alors au conoyau de la fl\`eche $K_R(X) \rightarrow K_R
^{-1}(X)$
qui s'ins\`ere dans la suite exacte de Bott (o\`u la derni\`ere fl\`eche est induite par la
r\'ealification des fibr\'es) :
$$K_R(X) \rightarrow  K_R
^{-1}(X) \rightarrow  K_C
^{-1}(X) \rightarrow K_R
^{1} (X)$$

\medskip
\noindent
L'obstruction \`a la transversalit\'e de $L_1$ et $L_2$ peut donc \^etre consid\'er\'ee comme un
\'el\'ement du groupe $N = Ker[K_C
^{-1}(X) \rightarrow K_R
^1 (X)]$, \'el\'ement qu'on notera $\gamma (E, L_1, L_2)$.

\medskip
\noindent
{\bf 7.} Remarque. Comme dans le paragraphe $3$, on peut montrer que le groupe $N \otimes Q $ s'identifie \`a
la somme des groupes de cohomologie de Cech $H^{4n-1}(X ; Q)$.

\medskip
\noindent
{\bf 8.} Soient $L_1$, $L_2$ et $ L_3$ trois lagrangiens tels que $L_1$ et $L_2$ d'une part, $L_2$ et $L_3$ d'autre
part, soient homotopiquement et stablement transversaux. Donc $L_3$ est d\'efini par un
morphisme de fibr\'es $ h : L = L_1\rightarrow L^* = L_2$ tel que $h^* = - h$. Pour que $L_3$ soit
homotopiquement transversal \`a $L_2$, il faut et il suffit donc que $ L$  puisse \^etre muni d'une
structure symplectique non d\'eg\'en\'er\'ee (ce qui est \'equivalent \`a une structure complexe).
Puisqu'on raisonne stablement, on peut remplacer $E$ par $E \oplus H(L)$, et donc $L$ par
$L\oplus L$. Sans restreindre la g\'en\'eralit\'e, on peut ainsi supposer que $ L$ est bien muni d'une
structure symplectique non d\'eg\'en\'er\'ee. On en d\'eduit le th\'eor\`eme suivant, dont la
d\'emonstration est calqu\'ee sur celle du th\'eor\`eme pr\'ec\'edent (en rempla\c  cant les m\'etriques
par des formes symplectiques non d\'eg\'en\'er\'ees).

\medskip
\noindent
\begin{theorem} Soit E un fibr\'e vectoriel muni d'une forme bilin\'eaire sym\'e\-trique non
d\'eg\'en\'er\'ee et soient $L_1$, ...,$ L_n $ des lagrangiens dans $E$. Pour que ces lagrangiens
soient homotopiquement et stablement transversaux, il faut et il suffit que $\gamma (E, L_1 , L_n)
= 0$ pour tout n (ce qui implique $\gamma(L_1, L_i, L_j) = 0$ pour tout couple $(i, j))$.
\end{theorem}

\medskip
\noindent
{\bf 9.} Remarque importante. Il est facile de voir que la th\'eorie $U(X)$ est un foncteur
semi-exact dont l'espace classifiant est le groupe orthogonal infini $ O = colim O(n)$
dans le cas $(+1)$-hermitien ou l'espace homog\`ene infini $U/O = colim U(n)/O(n)$ dans le
cas $(-1)$-hermitien. Si $X$ est un CW-complexe de dimension $N$, la transversalit\'e
homotopique stable est donc \'equivalente \`a la transversalit\'e homotopique si $ n > N+1$ dans le premier cas et $n > N+2$ dans le second (n \'etant le rang du fibr\'e $E$).

\medskip
\noindent
{\bf 10.} Revenons au cas d'un fibr\'e symplectique $E$ muni de trois lagrangiens $(L_1, L_2, L_3)$.
Sur le fibr\'e $L_1 \oplus L_2 \oplus L_3$ on peut  d\'efinir la forme quadratique suivante
$$q(x_1, x_2, x_3) = \phi(x_1, x_2) + \phi(x_2, x_3) + \phi(x_3, x_1)$$
utilis\'ee notamment par Leray \cite{Le} et Kashiwara \cite{S}.

\medskip
\noindent
\begin{theorem}. Les trois lagrangiens $L_1$, $L_2$ et $L_3$ sont transversaux si et seulement
si la forme quadratique q est non d\'eg\'en\'er\'ee.
\end{theorem}

\noindent
\begin{dem} Supposons d'abord que les trois lagrangiens ne soient pas
transversaux, par exemple que $L_1\cap L_2 \not =0$. Choisissons un \'el\'ement non nul $y_1 = y_2 \in
L_1\cap L_2$. Posons aussi $y_3 = 0$. Pour tout triplet $(x_1, x_2, x_3)$, calculons alors la forme
bilin\'eaire sym\'etrique associ\'ee \`a $q$, soit
$$\Psi(x, y) = \phi(x_1, y_2) - \phi(y_1, x_2) + \phi(x_2, y_3) - \phi(y_2, x_3) + \phi(x_3, y_1) - \phi(y_3, x_1)\;\; \;\;\;\; (F)$$

\medskip
\noindent
Elle se r\'eduit \`a $\phi(x_3, y_1) - \phi(x_3, y_2) + \phi(y_3, x_1) = \phi(x_3, y_1 - y_2) + \phi(y_3, x_1) = 0$. Ceci
montre que l'orthogonal de $F = L_1 \oplus L_2 \oplus L_3 $ pour la forme de Leray-Kashiwara
n'est pas nul et que la forme quadratique $q$ est d\'eg\'en\'er\'ee.

\medskip
\noindent
R\'eciproquement, supposons que les lagrangiens soient transversaux. On peut alors
poser $L_1= L, $ $ L_2 = L^*$, $ L_3$ \'etant d\'efini comme le graphe d'un morphisme $g : L \rightarrow L^*$
tel que $g^* = g.$ Consid\'erons le sous-espace $E_3$ de F form\'e des triplets $(x_1, x_2, x_3)$ tels
que $x_1 + x_2 - x_3 = 0$. Si on pose $x_3 = (u, gu)$ dans l'expression de $L_3$ comme sous-espace
de $L_1\oplus L_2$, on voit que $E_3$ est isomorphe \`a $L_1$. Pour calculer la forme $\Psi$
restreinte \`a $E_3$, on pose alors
$ x_1 = x'_1 = u,$
\noindent
$x_2 = x'_2 = gu,$
\noindent
$y_1 = y'_1= v,$
\noindent
 $y_2 = y'_2= gv.$

\medskip
\noindent
En appliquant la formule $(F)$ ci-dessus, on voit que la forme $\Psi$ restreinte \`a $E_3$ s'\'ecrit de
la mani\`ere suivante en fonction de $ u$ et $v \in L_1\cong E_3$, le symbole $< , >$ d\'esignant
l'accouplement entre un espace et son dual :

$$\Psi(u, v) = < u, gv > - <v, gu > + <v, gu> - < u, gv > + <v, gu > -$$
$$ < u, gv > = - 2 < u, gv >$$
Ainsi, au facteur $-2$ pr\`es, la forme bilin\'eaire sym\'etrique sur $E_3$ s'identifie \`a la forme
bilin\'eaire sym\'etrique associ\'ee \`a g (modulo l'identification entre $E_3$ et $L_1$). Par ailleurs,
si on consid\`ere de nouveau la formule $(F)$ ci-dessus, on voit que l'orthogonal de $H(L_1)
= L_1 \oplus L_2$ dans $F = L_1 \oplus L_2 \oplus L_3$ est form\'e des triplets  $(y_1, y_2, y_3)$, avec $ y_3 = ( y'_1
,  y'_2)$ tels que pour tout couple $(x_1, x_2)$, on ait $< x_1
, y_2 - y'_2) + <y_1 - y'_1,x_2> = 0$. Donc $y'_1
=
y_1$ et $ y'_2= y_2$, ce qui d\'etermine bien le sous-module $E_3$. Par cons\'equent, $F$ s'identifie \`a
la somme orthogonale de $ E_3$ et de $H(L_1)$ : ceci d\'emontre que la forme quadratique q est
non d\'eg\'en\'er\'ee.
\end{dem}

\medskip
\noindent
{\bf 11.} Consid\'erons maintenant le mono\" \i de des classes d'isomorphie des quadruplets
pr\'ec\'edents $(E, L_1, L_2, L_3)$ o\`u $L_1$, $L_2$ et $ L_3$ sont des lagrangiens transversaux. Notons $LK(X)$ le groupe de Grothendieck associ\'e \`a ce mono\" \i de (en r\'ef\'erence au travail de
Leray et Kashiwara).

\medskip
\noindent
\begin{theorem} Le groupe $LK(X)$ s'identifie au groupe de $K$-th\'eorie hermitienne de
l'espace $X$, c'est-\`a-dire \`a la somme directe de deux copies de la $K$-th\'eorie r\'eelle \footnote{ cf. par exemple le livre du premier auteur ``K-theory. An introduction'' p. $47$, exercice $9.22$ (Springer-
Verlag, $1978$)}.
De mani\`ere pr\'ecise, l'isomorphisme
$$LK(X) \rightarrow K_R(X) \oplus K_R(X)$$
est induit par la correspondance

$$(E, L_1, L_2, L_3) \mapsto [ \Psi] - H(L_1)$$
o\`u $\Psi $ est la forme quadratique de Leray-Kashiwara.
\end{theorem}

\noindent \begin{dem} D'apr\`es ce qui pr\'ec\`ede, les classes d'isomorphie des quadruplets
pr\'ec\'edents sont en correspondance bijective avec les formes bilin\'eaires sym\'etriques non
d\'eg\'en\'er\'ees $g : L_1\rightarrow L^*_1.$ Le th\'eor\`eme r\'esulte alors des consid\'erations du paragraphe $10$.

\end{dem}

\medskip
\noindent
{\bf 12.} G\'en\'eralisations diverses. Pour simplifier, nous nous sommes plac\'es dans le cadre g\'eo\-m\'etrique
des fibr\'es vectoriels. Les consid\'erations pr\'ec\'edentes se g\'en\'eralisent sans peine au cadre
des modules projectifs de type fini sur une $C^*$-alg\`ebre quelconque. Une autre
g\'en\'eralisation est de consid\'erer des fibr\'es \'equivariants  \footnote{ cf. l'article du second auteur ``Les invariants de Maslov \'equivariants" K-theory 25, p. 233-259 (2002).}. Les r\'esultats d\'emontr\'es
dans cette Note s'\'etendent aussi sans probl\`eme \`a ce cadre.

\def\refname{Bibliographie}

\end{document}